\def\lastrevised{June 7, 2011.}
\documentclass[12pt]{article}
\usepackage{amssymb,latexsym}

 \def\header{\today} 
\def\cutcrap{\addressarn\end{document}}

\def\al{\alpha}
\def\bb{{\mathcal B}}
\def\be{\beta}
\def\bspc{{\om^\om}}
\def\calc{{\mathcal C}}
\def\caln{{\mathcal N}}
\def\cc{{\mathcal  C}}
\def\cl(#1){{ cl}(#1)}
\def\concat(#1){\hat{\phantom{a}}\la #1\ra}
\def\cont{{\mathfrak c}}
\def\cspc{{2^\om}}
\def\de{\delta}
\def\dom(#1){dom(#1)}
\def\eps{\epsilon}

\def\gg{{\mathcal G}}
\def\hrusak{Hru\u{s}\'ak }

\def\isolated{\oplus 1}

\def\la{\langle}
\def\nicetree{{\rm NT}}
\def\om{\omega}
\def\pair(#1,#2){\la #1, #2 \ra}
\def\pf{\par\bigskip\noindent proof:\par}
\def\poset{{\mathbb P}}
\def\power{{\mathcal P}}

\def\proof{\pf}
\def\pr{\prime}
\def\qed{\par\noindent QED\par\bigskip}
\def\qq{{\mathbb Q}}
\def\ra{\rangle}
\def\res{{\upharpoonright}}
\def\rmand{\mbox{ and }}

\def\rmiff{\mbox{ iff }}

\def\root{{\la\ra}}
\def\rr{{\mathbb R}}
\def\si{\sigma}
\def\sm{\setminus}
\def\st{\;:\;} 
\def\su{\subseteq}
\def\symdiff{{\Delta}}
\def\uu{{\mathcal U}}
\def\zerovec{{\vec{0}}}

\newtheorem{theorem}{Theorem}
\newtheorem{lemma}[theorem]{Lemma}
\newtheorem{define}[theorem]{Definition}
\newtheorem{remark}[theorem]{Remark}
\newtheorem{prop}[theorem]{Proposition}
\newtheorem{ques}[theorem]{Question}
\newtheorem{cor}[theorem]{Corollary}
\newtheorem{example}[theorem]{Example}

\newenvironment{enum}
{\begin{list}{(\alph{enumi})}
{\usecounter{enumi}\setlength{\rightmargin}{\leftmargin}}}
{\end{list}}

\def\addressarn{\begin{flushleft}
Arnold W. Miller \\
miller@math.wisc.edu \\
http://www.math.wisc.edu/$\sim$miller\\
University of Wisconsin-Madison \\
Department of Mathematics, Van Vleck Hall \\
480 Lincoln Drive \\
Madison, Wisconsin 53706-1388 \\
\end{flushleft}}

\begin{document}

\begin{center}
{\large Uniquely Universal Sets}
\end{center}

\begin{flushright}
Arnold W. Miller\\
\end{flushright}

\begin{center}  Abstract\footnote{
\par Last revised \lastrevised}
\end{center}

\begin{quote}
We show that for any Polish space $Y$ that:

\centerline{ $\bspc\times Y$ has a uniquely universal set
iff $Y$ is not compact.}
\end{quote}

For any space $Y$ with a countable basis
there exists an open set
$U\su 2^\om\times Y$ which is universal for open subsets
of $Y$, i.e., $W\su Y$ is open iff there exists $x\in 2^\om$
with
$$U_x=^{{\rm def}}\{y\in Y\st (x,y)\in U\}=W.$$
To see this let $\{B_n\st n<\om\}$ be a basis for $Y$.
Define
$$(x,y)\in U \rmiff \exists n\; (x(n)=1 \rmand y\in B_n).$$

More generally if $X$ contains a homeomorphic copy of $\cspc$ then
$X\times Y$ will have a universal open set.

In 1995
Michael \hrusak mentioned the following problem to us.
Most of the results in this note were proved in
June and July of 2001.

\begin{center}
\hrusak's problem.
\end{center}

\begin{quote}
Let $X$, $Y$ be topological spaces, call $X$ the parameter space,
and $Y$ the base space.  When does there exists $U\subseteq X\times Y$
which is uniquely universal for the open subsets of $Y$? This
means the $U$ is open and for every open set $W\su Y$
there is a unique $x\in X$ such that $U_x=W$.
\end{quote}

Let us say that $X\times Y$ satisfies UU (uniquely universal property)
if there exists such an open set $U\su X\times Y$ which uniquely
parameterizes the open subsets of $Y$.  Note that the complement of
$U$ is a closed set which uniquely parameterizes the closed subsets
of $Y$.

\begin{prop}
(\hrusak) $2^\om\times 2^\om$ does not satisfy UU.
\end{prop}
\proof
The problem is the empty set.  Suppose $U$ is
uniquely universal for the closed subsets of $2^\om$.
Then there is an $x_0$ such that $U_{x_0}=\emptyset$
but all other cross sections are nonempty.  Take $x_n\to x_0$
but distinct from it.  Since all other cross sections are non-empty
we can choose $y_n\in U_{x_n}$.  But then  $y_n$ has a convergent
subsequence, say to $y_0$, but then $y_0\in U_{x_0}$.
\qed

More generally:

\begin{prop}\label{hrus2}
(\hrusak) Suppose $X\times Y$ has UU and $Y$ is compact.  Then
$X$ must have an isolated point.
\end{prop}
\pf
Suppose  $U\subseteq X\times Y$ witnesses
UU for closed subsets of $Y$ and $U_{x_0}=\emptyset$.
For every $y\in Y$ there exists
$U_y\times V_y$ open containing $(x_0,y)$ and missing $U$.
By compactness of $Y$ finitely many $V_y$ cover $Y$.
The intersection of the corresponding
$U_y$ isolate $x_0$.
\qed

Hence, for example,  $2^\om\times (\om+1)$,
$\bspc\times(\om+1)$, and $\bspc\times \cspc$ cannot
have UU.

\begin{prop}
Let $X=\cspc\isolated$ be obtained by attaching
an isolated point to $\cspc$.
Then $X\times \cspc$ has UU.
\end{prop}
\pf
Define $T\su 2^{<\om}$ to be a nice tree iff
\begin{enum}
\item $s\su t\in T$ implies $s\in T$ and
\item if $s\in T$, then either $s\concat(0)$
or $s\concat(1)$ in $T$.
\end{enum}
Let $\nicetree\su\power(2^{<\om})$ be the set of nice trees.
Define $U\su \nicetree\times \cspc$ by
$(T,x)\in U$ iff $\forall n \;\; x\res n\in T$.
Note that the empty tree $T$ is nice and parameterizes the empty set.
Also $\nicetree$ is a closed subset of $\power(2^{<\om})$ with
exactly one isolated point (the empty tree),
and hence it is homeomorphic to $X$.
\qed

\begin{ques}
Does $(\cspc \isolated)\times [0,1]$ have UU?
\end{ques}

\begin{remark}\label{powerset}
$2^\om\times \om$ has the UU property.  Just
let $(x,n)\in U$ iff $x(n)=1$.
\end{remark}

\begin{ques}
Does either $\rr\times \om $
or $[0,1]\times \om$ have UU?
Or more generally, is there any example of UU for a
connected parameter space?
\end{ques}

Recall that a topological space is Polish iff it is completely
metrizable and has a countable dense subset. A set is $G_\de$
iff it is the countable intersection of open sets. The
countable product of Polish spaces is Polish.  A $G_\de$ subset
of a Polish space is Polish (Alexandrov). 
A space is zero-dimensional iff it
has a basis of clopen sets. All compact zero-dimensional Polish
spaces without isolated points are homeomorphic to $\cspc$
(Brouwer). A zero-dimensional Polish space is homeomorphic to
$\bspc$ iff compact subsets have no interior
(Alexandrov-Urysohn).  For proofs of these facts see Kechris
\cite{kec} p.13-39.

\begin{define}
We use $\cl(X)$ to denote the closure of $X$.
\end{define}

\begin{prop} \label{locallycompact} Suppose $Y$ is Polish.
If $Y$ is locally compact but not compact,
then $\cspc\times Y$ has UU. So, for example,
$\cspc\times (\om\times 2^\om)$ and $\cspc\times \rr$
each have UU.
\end{prop}
\pf
Let $\bb$ be a countable base for $Y$ such that $b\in\bb$ implies
that $\cl(b)$ is compact.  Define $G\su\bb$ is good iff
$$G=\{b\in\bb \st \cl(b)\su \bigcup G\}$$
Let $\gg\su\power(\bb)$ be the family of good subsets of $\bb$.
We give the $\power(\bb)$ the topology from identifying
it with $2^{\bb}$.  Since $\bb$ is an infinite countable set
$\power(\bb)$ is homeomorphic to $\cspc$.  A sequence $G_n$
for $n<\om$ converges to
$G$ iff for each finite $F\su\bb$ we have
that $G_n\cap F=G\cap F$ for all but
finitely many $n$.   Hence $\gg$ is a closed
subset of $\power(\bb)$ since by compactness
$\cl(b)\su \bigcup G$
iff $\cl(b)\su \bigcup F$ for some finite $F\su G$.

There is a one-to-one correspondence between good families and
open subsets of $Y$:  Given any open set $U\su Y$ define
$$G_U=\{b\in\bb\st\cl(b)\su U\}$$
and given any good $G$ define $U_G=\bigcup G$.

We claim that no $G_0\in \gg$ is an isolated point.  Suppose
for contradiction it is.  Then there must be a basic open set $N$
with $\{G_0\}=\gg\cap N$.  A basic neighborhood has the following
form
$$N=N(F_0,F_1)=\{G\su \bb\st F_0\su G\rmand F_1\cap G=\emptyset\}$$
where $F_0, F_1\su \bb$ are finite.

For each $b\in F_1$ since
$G_0$ is good, $\cl(b)$ is not a subset of $\bigcup G_0$, and since
$\cl(\bigcup F_0)\su \bigcup G_0$, we can choose
a point $z_b\in \cl(b)\sm\cl(\bigcup F_0)$.
Since $Y$ is not compact, $Y\sm\cl(\bigcup (F_0\cup F_1))$ is nonempty.
Fix $z\in Y\sm\cl(\bigcup F)$.

Now let $U_1=Y\sm \{z_b\st b\in F_1\}$ and let $U_2=U_1\sm\{z\}$.
Then $G_{U_1},G_{U_2}$ are distinct elements of $N\cap\gg$.

Hence $\gg$ is a compact zero-dimensional metric space
without isolated points and therefore it is homeomorphic to $\cspc$.

To get a uniquely universal open set $U\su \gg\times Y$ define:
$$(G,x)\in U\rmiff \exists b\in G\;\; x\in b.$$

\qed

Is the converse true?  No.  Example \ref{examplepolish}
is a Polish space $Y$ such that $\cspc\times Y$ has
UU, but $Y$ is not locally compact.

\begin{lemma}\label{onto}
Suppose $f:X\to Y$ is one-to-one continuous and onto and
$Y\times Z$ has UU.  Then $X\times Z$ has UU.
\end{lemma}
\proof
Given $V\su Y\times Z$ witnessing UU,
let $$U=\{(x,y)\st (f(x),y)\in V\}.$$
\qed

Many uncountable standard Borel sets are
the bijective continuous image of the Baire space $\om^\om$.
According to the footnote on page 447 of Kuratowski \cite{kur}
Sierpinski proved in a 1929 paper that any
standard Borel set in which every point is a condensation point
is the bijective continuous image of $\bspc$.
We weren't able to find Sierpinski's paper but we
give a proof of his result in Lemma \ref{standardborel}.

We first need a special case for which we give a proof.

\begin{lemma}\label{bspcontocspc}
There is a continuous bijection $f:\bspc\to\cspc$.
\end{lemma}
\pf
Let $\pi:\om\to\om+1$ be a bijection.  It is automatically continuous.
It induces a continuous bijection $\pi:\bspc\to (\om+1)^\om$.
But $(\om+1)^\om$ is a compact Polish space without isolated points,
hence it is homeomorphic to $\cspc$.
\qed

Remark.  If $C\su\cspc\times\bspc$ is the graph of $f^{-1}$, then
$C$ is a closed set which uniquely parameterizes the family
of singletons of $\bspc$.

\begin{cor} \label{corparam}
If $\cspc\times Y$ has UU, then $\bspc\times Y$ has UU.
\end{cor}

\begin{ques}
Is the converse of Corollary \ref{corparam} false?
That is: Does there exist $Y$ such that $\bspc\times Y$ has UU but
$\cspc\times Y$ does not have UU?
\end{ques}

\begin{lemma}\label{bspconto}
Suppose $X$ is a zero-dimensional Polish space without isolated
points.  Then there exists a continuous bijection $f:\bspc\to X$.
\end{lemma}
\pf
Construct a subtree $T\su \om^{<\om}$ and  $(C_s\su X\st s\in T)$
nonempty clopen sets such that:
\begin{enumerate}
\item $C_{\root}=X$,
\item if $s\in T$ is a terminal node, then $C_s$ is
compact, and
\item if $s\in T$ is not terminal, then $s\concat(n)\in T$ for
every $n\in\om$ and $C_s$ is partitioned by
$(C_{s\concat(n)}\st n<\om)$
into nonempty clopen sets each of diameter less than
$\frac1{|s|+1}$.
\end{enumerate}
For each terminal node $s\in T$ choose a continuous bijection
$f_s:[s]\to C_s$ given by Lemma \ref{bspcontocspc}.
Define $f:\bspc\to X$ by
$f(x)=f_{x\res n}(x)$ if there exists $n$ such that
$x\res n$ is a terminal node of $T$ and otherwise determine
$f(x)$ by the formula:

$$\{f(x)\}=\bigcap_{n<\om} C_{x\res n}$$
Checking that $f$ is a continuous bijection is left to the
reader.
\qed

Remark. An easy modification of the above argument shows that any
zero-dimensional Polish space is homeomorphic to a closed
subspace of $\bspc$.  It also gives the classical result that
if no clopen sets are compact, then $X$ is homeomorphic to
$\bspc$.
A different proof of Lemma \ref{bspconto} is
given in Moschovakis \cite{mos} p. 12.

\begin{prop} \label{polish}
Suppose $Y$ is Polish but not compact.
Then $\bspc\times Y$ has the UU.
So for example, $\bspc\times\bspc$ and
$\bspc\times\rr$ both have UU.
\end{prop}
\pf
We assume that the
metric on $Y$ is complete and bounded.
Let $\bb$ be a countable basis for $Y$
of nonempty open sets which
has the property that no finite subset of $\bb$
covers $Y$.

For $s,t\in\bb$ define $t\lhd s$ iff $\cl(t)\subseteq s$
and diam$(t)\leq\frac12 $ diam$(s)$.

\begin{lemma}\label{lemprop}
Suppose $ G\subseteq \bb$ has the following properties:
\begin{enumerate}
  \item for all $t,s\in\bb$ if
  $t\subseteq s\in G$, then $t\in  G$ and
  \item $\forall s\in\bb\;$ if ($\forall t\lhd s\;\; t\in G$), then
$s\in  G$,
\end{enumerate}
then for any $s\in\bb$ if $s\subseteq\bigcup G$, then
$s\in G$.
\end{lemma}
\pf
Suppose not and $s\subseteq\bigcup G$ but $s\notin G$.
Note that there cannot be a sequence $(s_n:n\in\om)$
starting
with $s_0=s$, and with $s_{n+1}\lhd s_n$ and $s_n\notin  G$
for each $n$.  This is because
if $\{x\}=\bigcap_{n\in\om}s_n$, then $x\in s\subseteq\bigcup G$
and so for some $t\in G$ we have $x\in t$.  But then
for some sufficiently large $n$ we have that $s_n\subseteq t$
putting $s_n\in G$.   Hence there must be some
$t\unlhd s$ with $t\notin G$ but all $r\lhd t$ we have $r\in G$.
This is a contradiction.
\qed

Let $\gg\su\power(\bb)$ be the set of all $G$ which satisfy
the hypothesis of Lemma \ref{lemprop}.
Then we have a $G_\delta$ subset $\gg$ of $2^{\bb}$,
which uniquely parameterizes the open sets.
To finish the proof it is enough to see that no
point in $\gg$ is isolated.
If $\gg$ has an isolated point, there must be
$s_1,\ldots, s_n$ and $t_1,\ldots,t_m$ from $\bb$
such that
$$N=
\{G\in \gg\st s_1\in G,\ldots,s_n\in G,\;t_1\notin G,\ldots, t_m\notin G\}
$$
is a singleton.  Let $W=s_1\cup\cdots\cup s_n$. Since $N$ contains the
point $$G=\{s\in\bb\st s\su W\}$$
it must be that $N=\{G\}$.
For each $j$ choose $x_j\in t_j\sm W$.
Since the union of the $s_i$ and $t_j$ does not cover $Y$, we
can choose $$z\in Y\sm (s_1\cup\cdots \cup s_n\cup t_1\cup\cdots\cup t_m).$$
Take $r\in\bb$ with
$z\in r$ but $x_j\notin r$ for each $j=1,\ldots,m$.
Let $$G^\pr=\{s\in\bb\st s\su r\cup W\}.$$  Then $G^\pr\in N$ but
$G^\pr\neq G$ which shows that $N$ is not a singleton.
\qed

Proposition \ref{hrus2}  and \ref{polish} show that:

\begin{cor}
For $Y$ Polish:

\centerline{ $\bspc\times Y$ has UU iff $Y$ is not compact.}
\end{cor}

\begin{ques}\label{ques17}
Does $\cspc\times \bspc$ have UU?
\end{ques}

\begin{ques}
Does $\cspc\times \rr^\om$ have UU?
\end{ques}

Our next result follows from Proposition \ref{sigmacompact}
but has a simpler proof so we give it first.

\begin{prop}\label{ctbl}
If $X$ is a
countable metric space which is not compact,
then $\bspc\times X$ has UU. So, for example,
$\bspc\times\qq$ has UU.
\end{prop}

\proof
We produce a uniquely universal set for the open subsets
of $X$.

First note that there exists a countable basis $\bb$ for $X$
with the property that it is closed under finite unions and
$X\setminus B$ is infinite for every $B\in \bb$.  To see this
fix $\{x_n\st n<\om\}\su X$ an infinite set without a limit point,
i.e., an infinite closed discrete set.
Given a countable basis  $\bb$
replace it with finite unions of sets from
$$\{B\sm\{x_m:m>n\}\st n\in\om\rmand B\in\bb\}.$$
We may assume also that $\bb$ includes the empty set.

Next let
$$\poset=\{(B,F): B\in\bb,\;\;F\in[X]^{<\om},\rmand B\cap F=\emptyset\}.$$
Then $\poset$ is a partial order determined by
$p\leq q$ iff $B_q\su B_p$ and $F_q\su F_p$.
For $p\in\poset$ we write $p=(B_p,F_p)$. For $p,q\in\poset$ we
write $p\perp q$ to stand for $p$ and $q$ are incompatible, i.e.,
there does not exist $r\in \poset$ with $r\leq p$ and $r\leq q$.

We will code open subsets of $X$ by
good filters on $\poset$.
Define the family $\gg$ of good filters on $\poset$
to be the set of all $G\su\poset$ such that
\begin{enumerate}
\item $p\leq q$ and $p\in G$ implies $q\in G$,
\item $\forall p,q\in G$ exists $r\in G$ with $r\leq p$ and $r\leq q$,
\item $\forall x\in X\;\exists p\in G \;\; x\in B_p\cup F_p$, and
\item $\forall p\in\poset$ either $p\in G$ or $\exists q\in G\;\; p\perp q$.
\end{enumerate}
Since the poset $\poset$ is countable we can identify
$\power(\poset)$ with $\power(\om)$ and hence $\cspc$.
We give $\gg\su\power(\poset)$
the subspace topology.  Note that $\gg$ is $G_\de$ in this
topology.  Note also that the sets
$$[p]=\{G\in\gg\st p\in G\}$$
form a basis for $\gg$ (use conditions (2) and (4)
to deal with finitely many $p_i$ in $G$ and finitely
many $q_j$ not in $G$).

Note that since $X\sm B_p$ is always infinite, for any $p\in\poset$
there exists $r,q\leq p$ with $r\perp q$.  Namely, for some
$x\in X\sm (B_p\cup F_p)$ put $x$ into $B_r\cap F_q$.
It follows that no element of $\gg$ is isolated.
So $\gg$ is a zero-dimensional Polish space without
isolated points.
Hence by Lemma \ref{bspconto} there is a continuous
bijection $f:\bspc\to\gg$.

For $G\in \gg$, let
$$U_G=\bigcup \{B_p\st p\in G\}.$$
For any $U\su X$ open, define
$$G_U=\{p\in\poset\st B_p\su U \rmand F_p\cap U=\emptyset\}.$$
The maps $G\to U_G$ and $U\to G_U$ show that
there is a one-to-one correspondence between $\gg$ and the
open subsets of $X$.

Finally define $\uu\su \gg\times X$ by
$$(G,x)\in \uu \rmiff \exists p\in G \;\; x\in B_p.$$
This witnesses the UU property for $\gg\times X$ and
so by Lemma \ref{onto}, we have UU for $\bspc\times X$.

\qed

\begin{ques}
Does $\cspc\times\qq$ have UU?
\end{ques}

Our next result Proposition \ref{sigmacompact} implies Proposition
\ref{ctbl} but needs the following Lemma:

\begin{lemma}\label{standardborel}(Sierpinski)
Suppose $B$ is a Borel set in a Polish space for which every point is
a condensation point.  Then there exists a continuous
bijection from $\bspc$ to $B$.
\end{lemma}
\pf
We use that every Borel set is the bijective image of a closed
subset of $\bspc$.  This is due to Lusin-Souslin
see Kechris \cite{kec} p.83 or Kuratowski-Mostowski \cite{kurmos}
p.426.

Using the fact that every uncountable Borel set contains a perfect
subset it is easy to
construct $K_n$ for $n<\om$ satisfying:
\begin{enumerate}
\item $K_n\su B$ are pairwise disjoint,
\item $K_n$ are homeomorphic copies of $\cspc$ which are
nowhere dense in $B$, and
\item every nonempty open subset of $B$ contains
infinitely many $K_n$.
\end{enumerate}

Let $B_0=B\sm\bigcup_{n<\om}K_n$.  We may assume $B_0$ is
nonempty, otherwise just split $K_0$ into
two pieces.  Since it is a Borel set,
there exists $C\su\bspc$ closed and a continuous bijection
$f:C\to B_0$.  Define
$$\Gamma=\{s\in \om^{<\om}\st [s]\cap C=\emptyset \rmand [s^*]\cap C\neq
\emptyset\}$$
where $s^*$ is the unique $t\su s$ with $|t|=|s|-1$.  Without loss
we may assume that $C$ is nowhere dense and hence $\Gamma$ infinite.
Let $\Gamma=\{s_n\st n<\om\}$ be a one-one enumeration.
Note that $\{C\}\cup\{[s_n]\st n<\om\}$ is a partition of $\bspc$.

Inductively choose
$l_n>l_{n-1}$ with $K_{l_n}$ a subset of the
ball of radius $\frac1{n+1}$ around $f(x_n)$ for some
$x_n\in C\cap [s_n^*]$.
For each $n<\om$ let
$f_n:[s_n]\to K_{l_n}$ be a continuous bijection.

\bigskip

Then $g=f\cup\bigcup_{n<\om}f_n$ is a continuous bijection
from $\bspc$ to $B_0\cup\bigcup_{n<\om}{K_{l_n}}$.

\bigskip\noindent
To see that it is continuous suppose for contradiction
that $u_n\to u$ as $n\to\infty$ and $|g(u_n)-g(u)|>\eps>0$
all $n$.  Since $C$ is closed if infinitely many $u_n$ are
in $C$, so is $u$ and we contradict continuity of $f$.  If
$u\in [s_n]$, then we contradict the continuity of $f_n$,
So, we may assume that all $u_n$ are not in $C$ but $u$ is in
$C$.  By the continuity of $f$ we may find $s\su u$ with
$f([s]\cap C)$ inside a ball of radius $\frac\eps{3}$ around $f(u)$.
Find $n$ with $\frac1{n+1}<\frac\eps3$ for which there is
$m$ such that $u_m\in [s_n]$ and $s_n\supseteq s$.
But then $g(u_m)=f_n(u_m)\in K_{l_n}$ and
$K_{l_n}$ is in a ball of radius $\frac1{n+1}$ around
some $f(x_n)$ with $x_n\in [s_n^*]\cap C$.
This is a contradiction:
$$d(g(u),g(u_m))\leq d(f(u),f(x_n))+d(f(x_n),f_n(u_m))\leq \frac23\eps.$$

Next let $I=\om\sm\{l_n\st n<\om\}$.  Then there exists continuous
bijection $$h:I\times\bspc\to \bigcup_{i\in I}K_i.$$
Finally
$g\cup h$ is a continuous bijection from $\bspc \oplus (I\times\bspc)$
to $B_0\cup \bigcup_{n<\om}K_n=B$.  Since $\bspc \oplus (I\times\bspc)$
is a homeomorphic copy of $\bspc$ we are done.
\qed

\begin{prop} \label{sigmacompact}
$\bspc\times Y$ has UU for any $\si$-compact but not
compact subspace $Y$ of a Polish
space.  So for example, $\bspc\times (\qq\times \cspc)$ has UU.
\end{prop}
\pf
Let $Y=\bigcup_{n<\om}K_n$ where each $K_n$ is compact.  Since
$Y$ is not compact we can choose a countable basis $\bb$ for
$Y$ such that for any finite $F\su\bb$ the closed
set $Y\sm\bigcup F$ is not compact.
Define $G\su\bb$ to be good iff for
every $b\in\bb$ if $\cl(b)\su\bigcup G$ then $b\in G$.
Let $\gg\su\power(\bb)$ be the family of good sets.

Note that $\gg$ is a ${\bf \Pi}^0_3$ set:
$$G\in\gg \rmiff
\forall b\in\bb \;\;(\forall n \;\;\cl(b)\cap K_n\su \bigcup G)\to b\in G$$
Note that $\cl(b)\cap K_n\su\bigcup G$ iff there is a finite
$F\su G$ with $\cl(b)\cap K_n\su\bigcup F$.

To finish the proof it is necessary
to see that basic open sets
in $\gg$ are uncountable.
Given $b_i,c_j\in\bb$ for $i<n$ and $j<m$ suppose that
$$N=
\{G\in \gg\st b_0\in G,\ldots,b_{n-1}\in G,\;
\;c_0\notin G,\ldots, c_{m-1}\notin G\}
$$
is nonempty.
Since it is nonempty we
can choose points $$u_j\in \cl(c_j)\sm \bigcup_{i<n}\cl(b_i)$$
for $j<m$.
Since
$$Y\sm \cl(\bigcup_{i<n,j<m}b_i\cup c_j)$$
is not compact, it contains an infinite discrete closed
set $Z=\{z_n\st n<\om\}$.  But then given any $Q\su Z$ we can
find an open set $U_Q$ with $\bigcup_{i<n}\cl(b_i)\su U_Q$,
$u_j\notin U_Q$ for $j<m$,
and $U_Q\cap Z=Q$.  Let
$$G_Q=\{b\in\bb\st \cl(b)\su U_Q\}$$
Since each $G_Q\in N$ we have that $N$ is uncountable.
By Lemma \ref{standardborel} and \ref{onto}, we are done.
\qed

\begin{ques}\label{fsigma}
For what Borel spaces $Y$ does $\bspc\times Y$ have UU?
\par For example, does $\bspc\times \qq^\om$ have UU?
\end{ques}

\begin{prop}
For every $Y$ a $\bf\Sigma^1_1$ set there exists a $\bf\Sigma^1_1$
set $X$ such that $X\times Y$ has UU.
\end{prop}
\pf
Suppose $Y\su Z$ where $Z$ is Polish and $\bb$ is a countable
base for $Z$.  Define $\gg\su\power(\bb)$ by
$$G\in \gg \rmiff \forall b\in\bb \;\;[(b\cap Y\su\bigcup G)\to b\in G]$$
Note that $b\cap Y\su\bigcup G$ is $\bf\Pi^1_1$
and so $\gg$ is $\bf\Sigma^1_1$.
\qed

We use the next Lemma for Example \ref{examplepolish}.

\begin{lemma} \label{omegatimescspc}
For any space $Y$

$(\om\times\cspc)\times Y$ has UU $\;\;\;$ iff $\;\;\cspc\times Y$ has UU.
\end{lemma}
\proof
Suppose $C\su (\om\times\cspc)\times Y$ is a closed set uniquely
universal for the closed subsets of $Y$.

Since the whole space $Y$ occurs as
a cross section of $C$ without loss we may assume
that $Y=C_{(0,\zerovec)}$ where $\zerovec$ is the constant zero
function.

For each $n>0$ let
$$K_n=\{(0^{n}
\concat(\star,x_0,x_1,\ldots),y)\in\cspc\times Y\st ((n,x),y)\in C\}$$
By $0^{n}\concat(\star,x_0,x_1,\ldots)$ 
we mean a sequence of $n$ zeros followed by the
special symbol $\star$ and then the (binary) digits of $x$.
Note that the $K_n$ converge to
$\zerovec$.
Let
$$K_0=\{(x,y)\st ((0,x),y)\in C\}$$
Let $B=\bigcup_{n<\om}K_n\su S\times Y$ where
$$S=\cspc\cup\bigcup_{n>0}\{0^{n}\concat(\star,x_0,x_1,\ldots)
\st x\in\cspc\}.$$
Note that $S$ is homeomorphic to $\cspc$ and there is a one-to-one
correspondence between the cross sections of $B$ and
cross sections of $C$.
Note that $B$ is closed in $S\times Y$:
If $(x_n,y_n)\in B$ for $n<\om$ is a sequence
converging to $(x,y)\in S\times Y$ and $x$ is not the zero vector it
is easy to see that $(x,y)\in B$.  On the other hand if $x$
is the zero vector, then since
$B_\zerovec=Y$, it is automatically
true that $(x,y)\in B$.

Hence UU holds for $\cspc\times Y$.

For the opposite direction just note that
$(\om+1)\times \cspc$ is homeomorphic
to $\cspc$ and
there is a continuous bijection from $\om\times\cspc$ onto
$(\om+1)\times \cspc$.

\qed

Next we describe our counterexample to the converse of
Proposition \ref{locallycompact}. 
Let $Z=(\om\times\om)\cup \{\infty\}$.
Let each $D_n=\{n\}\times\om$ be an
infinite closed discrete
set and let the sequence of $D_n$ ``converge'' to $\infty$, 
i.e., each neighborhood
of $\infty$ contains all but finitely many $D_n$.
Equivalently $Z$ is homeomorphic to:
$$Z^\pr=\{x\in\bspc \st |\{n\st x(n)\neq 0\}|\leq 1\}.$$
The point $\infty$ is the constant zero map, while
$D_n$ are the points in $Z^\pr$ with $x(n) \neq 0$.
Note that $Z^\pr$ is a closed subset of
$\bspc$, hence it is Polish.
This seems to be the simplest nonlocally compact Polish space.

\begin{example} \label{examplepolish}
$Z$ is a nonlocally compact Polish space such that
$2^\om\times Z$ has the UU.
\end{example}
\pf
$Z=\bigcup_{n<\om}D_n\cup\{\infty\}$.  Note
that $X\su Z$ is closed iff $\infty\in X$ or
$X\su\bigcup_{i\leq k}D_i$ for some $k<\om$.
By Lemma \ref{omegatimescspc} it is enough to see that
$(\om\times 2^\om)\times Z$ has the UU.

Let $P_n=\{n\}\times 2^\om$.

Use $P_0$ to uniquely parameterize all
subsets of $Z$ which contain the point at infinity,
see Remark \ref{powerset}.

Use $P_1$ to uniquely parameterize all $X\su D_0$,
including the empty set.

For $n=1+2^{k-1}(2l-1)$ with $k,l>0$ use $P_n$ to uniquely 
parameterize all
$X\su \bigcup_{i\leq k}D_i$ such that $D_k$ meets $X$ and
the minimal element of $D_k\cap X$ is the $l$-th element
of $D_k$.

\qed

Our next two results have to do with Question \ref{ques17}.

\bigskip
\begin{prop}
Existence of UU for $\cspc\times \bspc$ is equivalent to:
\par There exist a $\cc\su \power(\om^{<\om})$ homeomorphic
to $\cspc$
such that every $T\in \cc$ is a subtree of $\om^{<\om}$
(possibly with terminal nodes) and such that for every closed
$C\su\bspc$ there exists a unique
$T\in\cc$ with $C=[T]$.
\end{prop}
\pf
Given $C\su \cspc\times\bspc$ witnessing UU
for closed subsets of $\bspc$.
Let
$$[T]=\{(s,t):([s]\times [t])\cap C \neq \emptyset\}.$$
Define $f:\cspc\to \power(\om^{<\om})=2^{\om^{<\om}}$ by
$f(x)(s)=1$ iff $(x\res n,s)\in T$ where $n=|s|$.
Then $f$ is continuous, since $f(x)(s)$ depends
only on $x\res n$ where $n=|s|$.

If $T_x=f(x)$, then $[T_x]=C_x$.  Hence $f$ is one-to-one,
so its image $\cc$ is as described.

\qed
\bigskip

\begin{prop}
Suppose $\bspc\times Y$ has UU.  Then there exists $U\su \cspc\times Y$
an $F_\si$ set such that all cross
sections $U_x$ are open and for every open $W\su Y$ there
is a unique $x\in \cspc$ with $U_x=W$.
\end{prop}
\pf
Let $\om\isolated$ denote the discrete space with one isolated
point adjoined and let $\om+1$ denote the compact space consisting
of a single convergent sequence.  Then $(\om\isolated)^\om$ is
homeomorphic to $\bspc$ and $(\om+1)^\om$ is homeomorphic to $\cspc$.

Assume that $U\su (\om\isolated)^\om \times Y$ is an
open set witnessing UU.  Then $U$ is an $F_\si$ set in
$(\om+1)^\om\times Y$.   To see this suppose that
$U=\bigcup_{n<\om} [s_n]\times W_n$.
Note that $[s_n]$ is closed in $(\om+1)^\om$.
\qed

\bigskip

Here $\oplus$ refers to the topological sum of
disjoint copies or equivalently a clopen separated union.
\begin{prop}
Suppose $X_i\times Y_i$ has UU for $i\in I$.
Then
$$(\prod_{i\in I}X_i)\times (\bigoplus_{i\in I} Y_i)\mbox{ has UU.}$$
So, for example, if $\cspc\times Y$ has UU, then
$\cspc\times (\om\times Y)$ has UU.
\end{prop}
\pf
Define
$$((x_i)_{i\in I},y)\in U \rmiff \exists i\in I\;\; (x_i,y)\in U_i$$
where the $U_i\su X_i\times Y_i$ witness UU.
\qed
\bigskip

Except for Proposition \ref{hrus2} we have given no negative results.
The following two propositions are the best we could do in
that direction.

\begin{prop}
Suppose $U\su X\times Y$ is an open set universal for the
open subsets of $Y$.  If $X$ is second countable, then so
is $Y$.
\end{prop}

\pf
$U$ is the union of open rectangles of the form $B\times C$
with $B$ open in $X$ and $C$ open in $Y$.  Clearly
we may assume that $B$ is from a fixed countable basis for
$X$.  Since $\bigcup_i B\times C_i=B\times \bigcup_i  C_i$
we may write $U$ as a countable union:
$$U=\bigcup_{n<\om} B_n\times C_n$$
where the $B_n$ are basic open sets in $X$ and the
$C_n$ are open subsets of $Y$.   But this implies
that $\{C_n\st n<\om\}$ is a basis for $Y$ since
for each $x\in X$
$$U_x=\bigcup \{C_n\st x\in B_n\}$$
\qed

\begin{prop}
There exists a partition $X\cup Y=\cspc$ into Bernstein sets
$X$ and $Y$ such that for every Polish space $Z$ neither
$Z\times X$ nor $Z\times Y$ has UU.
\end{prop}

\pf
Note that up to homeomorphism there are only continuum
many Polish spaces.
If there is a UU set for $Z\times X$, then there is
an open $U\su Z\times\cspc$ such that $U\cap (Z\times X)$
is UU.  Note that the cross sections of $U$ 
must be distinct open subsets of $\cspc$.
Hence it suffices to prove the
following:

Given $\uu_\al$ for $\al<\cont$
such that each $\uu_\al$ is a family of open subsets
of $\cspc$ either
\begin{enum}
\item there exists distinct $U,V\in\uu_\al$ with
$U\cap X=V\cap X$ or
\item there exists $U\su\cspc$ open such that $U\cap X\neq V\cap X$
for all $V\in\uu_\al$.
\end{enum}
And the same for $Y$ in place of $X$.

Let $P_\al$ for $\al<\cont$ list all perfect subsets of
$\cspc$ and  let $\{z_\al\st\al<\cont\}=\cspc$.
Construct $X_\al,Y_\al\su\cspc$ with
\begin{enumerate}
\item $X_\al\cap Y_\al=\emptyset$
\item $\al<\be$ implies $X_\al\su X_\be$ and $Y_\al\su Y_\be$
\item $|X_\al\cup Y_\al|=|\al|+\om$
\item If there exists distinct $U,V\in\uu_\al$ such
that $U\symdiff V$ is a countable set disjoint
from $X_{\al}$, then there exists such a pair with
$U\symdiff V\su Y_{\al+1}$
\item
If there exists distinct $U,V\in\uu_\al$ such
that $U\symdiff V$ is a countable set disjoint
from $Y_{\al+1}$, then there exists such a pair with
$U\symdiff V\su X_{\al+1}$
\item $P_\al$ meets both $X_{\al+1}$ and $Y_{\al+1}$
\item $z_\al\in (X_{\al+1}\cup Y_{\al+1})$
\end{enumerate}
First we do (4) then (5) and then take care of (6) and (7).

Let $X=\bigcup_{\al<\cont}X_\al$ and
$Y=\bigcup_{\al<\cont}Y_\al$.

Fix $\al$ and let us verify (a) or (b) holds.
Take any point $p\in X\sm X_{\al+1}$.  If (b) fails there
must be $U,V\in\uu_\al$ with $X=X\cap U$ and
$X\sm\{p\}=X\cap V$.  Then $(U\symdiff V)\cap X_{\al+1}=\emptyset$.
Since $X$ is Bernstein and $(U\symdiff V)\cap X$ has only
one point in it, it must be that $U\symdiff V$ is countable.
Then by our construction we have chosen distinct $U,V\in\uu_\al$
with $U\symdiff V\su Y$ therefor $U\cap X=V\cap X$, so (a) holds.

A similar argument goes through for $Y$.

\qed

Finally, and conveniently close to the bibliography, we note some
papers in the literature which are related to the property UU.
Friedberg \cite{friedberg} proved that there is one-to-one
recursively enumerable listing of the recursively enumerable sets.
This is the same as saying that there is a (light-face)
$\Sigma_1^0$ subset $U\su\om\times\om$ which is uniquely
universal for the $\Sigma_1^0$ subsets of $\om$.
Brodhead and Cenzer \cite{brodhead} prove the analogous
result for (light-face) $\Sigma_1^0$ subsets of $\cspc$.

Becker \cite{becker} considers unique parameterizations of
the family of countable sets by Borel or
analytic sets.  Gao, Jackson, Laczkovich, and Mauldin
\cite{gao} consider several other problems of unique parameterization.

\cutcrap

\newpage

\begin{center}
{\bf Miscellaneous Crapolo}
\end{center}

How about
Unique Parameterization (UP) for other families of sets?
It was remarked
that if $f:\bspc\to\cspc$ is 1-1,onto, continuous,
then
 $$U=\{(x,y): f(y)=x\}\subseteq \cspc\times\bspc$$
is UP for the singleton subsets of $\bspc$.

How about doubleton s? Does there exist a closed set
$C\subseteq\cspc\times\bspc$
such that $|C_x|=2$ for every $x$ and for every $F\subseteq\bspc$
of size 2 there exists an $x$ with $C_x=F$?  (If so,
can we get it to be UP?).

\bigskip
Note that we can do singletons (uniquely) and using that we can do
$F$ of size 1 or 2 uniquely.  For $F$ of size 3 can do but I don't
know about uniquely.  Can also use
$\om\times\cspc$ or $\bspc$ as unique parameters for doubleton s.

\bigskip
What about other classes of closed sets.  For example, the perfect
subsets of $\cspc$?  Or the family of finite closed sets.
Or the family of infinite discrete subsets of $\bspc$?

\bigskip
\bigskip
For $X$ a Polish space let $C(X)$ be the compact non-empty subsets
of $X$ with the Hausdorff metric:
$$d(K,L)=sup_{y\in L}dist(y,K)\;\;+\;\;sup_{x\in K}dist(x,L)$$
I think this is always a Polish space.  If we define
$Q\su C(X)\times X$ by $(K,x)\in Q$ iff $x\in K$ we get
a uniquely universal set for the compact subsets of $X$.

\bigskip
Is being a uniquely universal or universal set open sets
absolute for Polish products $U\su X\times Y$?

\bigskip

Can $U\subseteq X\times Y$ witness simultaneously
the UU for each?  Or even just universal without
the unique?  No it cannot.  If $U_x=Y$ then $x$ is in every $U^y$.

\bigskip
Investigate UU for maps $f:Y\to Z$, i.e.
$F:X\times Y\to Z$ giving every map as a unique cross section.

\bigskip
What about doubly universal uniqueness sets?  Suppose
$U\su X\times Y$ is a uniquely universal set.
Define $V,W\su (X\times X)\times Y$ by
$$V=\{((x,y),z): (x,z)\in U\}\rmand W=\{((x,y),z): (y,z)\in U\}$$
Then for any pair of open sets $A,B\su Y$ there is a unique
$z=(x,y)$ with $V_z=A$ and $W_z=B$.

\bigskip
On question about $\cspc\times\bspc$: what kind of function
is $T\to \hat{T}$ where $\hat{T}=\{s\in T: T_s$ has branch$\}$?

\bigskip
Its $\Sigma_1^1$.  Note that it is one-to-one onto
from a putative universal set to Trees without
terminal nodes (except root).

\bigskip\bigskip
On the question:
Does there exists a closed (perfect?)
$\calc\subseteq P(\om^{<\om})$ of trees
such that for every $C\subseteq\bspc$ closed there exists a unique
$T\in\calc$ with $[T]=C$?

Various thought on the matter:
$$\hat{T}=\{s\in T:[s]\cap[T]\not=\emptyset\}$$
$$D(T)=\{s\in T:\exists n \; sn\in T\}$$

$\caln$ subtrees of $\om^{<\om}$ with no terminal nodes
equivalently fixed points of $D$.

Note that $\calc\to\caln$ via $T\to \hat{T}$ is 1-1 and onto.

$\caln\to cl(\bspc)$ via $T\to[T]$ is one-one and onto.

\addressarn

\end{document}